\documentclass[oneside,a4paper,11pt,notitlepage]{article}
\usepackage[left=0.8in, right=0.8in, top=1.5in, bottom=1.5in]{geometry}

\usepackage[T1]{fontenc} 
\usepackage[utf8]{inputenc} 
\usepackage[english]{babel}
\usepackage{amssymb}
\usepackage{lipsum} 
\usepackage{lmodern}
\usepackage{amsmath}
\usepackage{amsthm}
\usepackage{bm}
\usepackage{mathtools}
\usepackage{braket}
\usepackage{esint}
\usepackage{comment}
\newcommand{\abs}[1]{{\left|#1\right|}}
\newcommand{\norma}[1]{{\left\Vert#1\right\Vert}}

\usepackage{relsize}

\usepackage{enumitem}
\usepackage{booktabs}
\usepackage{graphicx}
\usepackage{tikz}
\usetikzlibrary{patterns}
\usepackage{multicol}
\usepackage{caption}
\usepackage{bigints}
\usepackage[skins,theorems]{tcolorbox}
\tcbset{highlight math style={enhanced,
colframe=black,colback=white,arc=0pt,boxrule=1pt}}
\def\XXint#1#2#3{{\setbox0=\hbox{$#1{#2#3}{\int}$}
    \vcenter{\hbox{$#2#3$}}\kern-.5\wd0}}

\theoremstyle{definition}
\newtheorem{definizione}{Definition}[section]
\theoremstyle{plain}
\newtheorem{teorema}{Theorem}[section]

\newtheorem{lemma}[teorema]{Lemma}
\newtheorem{prop}[teorema]{Proposition}
\newtheorem{corollario}[teorema]{Corollary}

\theoremstyle{definition}
\newtheorem{esempio}{Example}[section]
\newtheorem{oss}[esempio]{Remark}

\newtheorem*{open*}{Open problems}
\DeclareMathOperator{\R}{\mathbb{R}}

\makeatletter
\newcommand{\myfootnote}[2]{\begingroup
	\def\@makefnmark{}%
	\addtocounter{footnote}{-1}%
	\footnote{\textbf{#1} #2}
	\endgroup}
\makeatother

\definecolor{mygreen}{RGB}{0, 255, 0}

\usepackage{hyperref}
\hypersetup{linktoc=none, bookmarksnumbered, colorlinks=true, linkcolor=magenta, citecolor=cyan}
\title{A quantitative Talenti comparison principle for the Monge-Ampère equation}
\author{Antonio Porricelli*, Francesco Salerno}
\date{}

\newcommand{\Addresses}{{
\bigskip 
  
   \medskip

    \textit{E-mail address}, A.~Porricelli* (corresponding author) \texttt{antonio.porricelli2@unina.it} 
  
   \medskip 
 \textsc{Università degli studi di Napoli "Federico II", Dipartimento di Matematica e Applicazioni "R. Caccioppoli", Via Cintia, Complesso Universitario Monte S.
Angelo, 80126 Napoli, Italy.}
   
   \medskip 
     \textit{E-mail address}, F.~Salerno: \texttt{f.salerno@ssmeridionale.it} 
   \medskip 
   
 \textsc{Mathematical and Physical Sciences for Advanced Materials and Technologies, Scuola Superiore Meridionale, Via Mezzocannone 4, 80134 Napoli, Italy.}

 \par\nopagebreak 

}} 

\begin{document}
\maketitle
\begin{abstract}
\noindent We prove a quantitative version of the comparison principle for the Monge-Ampère equation established by Talenti in dimension two and generalized by Tso to $k$-Hessian equations in arbitrary dimension. Namely, we estimate the $L^\infty$ distance between the $(n-1)$-symmetrization of the solution and the solution of the symmetrized problem in terms of the Hausdorff asymmetry of the domain. This extends to the case $k=n$ a recent result for $k<n$. As a consequence, we obtain a quantitative version of a sharp ABP-type estimate in terms of the diameter of the domain.
\end{abstract}
\section{Introduction}
In recent decades, symmetrization techniques have been extensively studied in relation to the qualitative properties of solutions to second-order elliptic boundary value problems. In the seminal paper \cite{talenti76}, Talenti uses the Schwarz symmetrization to prove sharp a priori bounds for solutions to the Poisson equation with Dirichlet boundary conditions (see \cite{Alvino_Lions_Trombetti_1986} for the rigidity of this result). A key feature of Schwarz symmetrization is that it preserves the volumes of the superlevel sets of a function, producing a radially symmetric function (see \cite{Kawohl, Kesavan} for further details on symmetrizations). In \cite{talenti1981}, the author introduces a new type of symmetrization that preserves the perimeter of the level sets of a function in the plane. This idea has been extended in \cite{tso1989}, where K. Tso introduces a symmetrization that preserves the $k$-th quermassintegral of the level sets of a function, producing a radially symmetric function. The corresponding result is stated below.\\
Let $\Omega\subset\mathbb{R}^n$, $n\ge2$, be an open, bounded and convex set, and let $f:\Omega\times[0,+\infty)\to(0,+\infty)$ be measurable, nonincreasing with
respect to its second variable, and such that $f_0:=f(\cdot,0)\in L^1(\Omega)$.
Let $u\in C^2(\Omega)\cap C(\overline\Omega)$ be a concave solution to
\begin{equation}\label{MAequation}
\begin{cases}
\det(-D^2u)=f(x,u) & \text{in }\Omega,\\
u=0 & \text{on }\partial\Omega .
\end{cases}
\end{equation}
Let $\Omega^*_{n-1}$ be the ball centered at the origin with the
same $(n-1)$-quermassintegral as $\Omega$, let $u^\# _{n-1}$ be the
$(n-1)$-symmetrization of $u$ (see Section \ref{Sec2} for the definitions), and let $v$ be the radial solution to
\begin{equation}\label{MAequationSymm}
\begin{cases}
\det(-D^2v)=f^\#_0 & \text{in }\Omega^*_{n-1},\\
v=0 & \text{on }\partial\Omega^*_{n-1},
\end{cases}
\end{equation}
where $f^\#_0$ denotes the Schwarz rearrangement of $f_0$, extended by zero outside the ball having the same measure as $\Omega$.
The following comparison principle was proved in \cite{talenti1981} for $n=2$ and extended in \cite{tso1989} for $n\geq2$:
\begin{equation}\label{Tso}
    u^\# _{n-1}(x) \leq v(x), \quad \forall x \in \Omega^*_{n-1}.
\end{equation}
Moreover, equality holds in \eqref{Tso} if and only if there exists $x_0\in\R^n$ such that $\Omega=\Omega^*_{n-1}+x_0$, $u(x)=u^\#_{n-1}(x -x_0)$ for every $x\in\Omega$ and $f(x,u(x))=f_0(x)=f_0^\#(x-x_0)$ for almost every $x\in\Omega$. Actually the inequality \eqref{Tso} is the special case of a family of inequalities proved in \cite{tso1989}. More precisely, Tso proved the inequality for a more general differential operator, namely the $k$-Hessian operator, that interpolates the Laplace operator ($k=1$) and the Monge-Ampère operator ($k=n$) (see \cite{MS2} for more details on the $k$-Hessian operator).\\
In recent years, the study of quantitative improvements to inequalities of the form \eqref{Tso} has attracted the attention of many mathematicians. In \cite{ABMP} the authors studied the case of the Laplace operator (see also \cite{AcamporaLamboley, AcamporaLmaboleyMazari} for an improvement of the quantitative inequality with an optimal exponent of the asymmetry).
In \cite{MS} a quantitative version of \eqref{Tso} for the $k$-Hessian operator, with $k<n$, has been proved. The idea behind their proof is an application of a quantitative version of a geometric inequality, namely the Alexandrov-Fenchel inequality. However, this method fails in the case $k=n$, in which this geometric inequality is no longer applied to the level sets of the function. In this work, we extend the result to the case $k=n$. The idea that allows us to overcome this difficulty is as follows: we obtain the quantitative result we are interested in by estimating H\"older’s deficit arising in the original proofs of Talenti and Tso. We then estimate this remainder in terms of a suitable geometric asymmetry, that is, the Hausdorff asymmetry (see Definition \ref{Hauss}) $\alpha_\mathcal{H}(\Omega)$, an index that measures the $L^\infty$ distance of $\Omega$ to a ball. More precisely, we prove the following.
\begin{teorema}\label{Main Teo}
    Let $\Omega\subset\mathbb{R}^n$, $n\ge2$, be an open, bounded and convex set. Let $u$, $v$ be the solutions to \eqref{MAequation}, \eqref{MAequationSymm} respectively. Then
    \begin{equation}\label{RisFinale}
        \frac{\|v-u_{n-1}^\# \|_{L^{\infty}(\Omega^*_{n-1})}}{\|u\|_{L^\infty(\Omega)}}\geq\frac{\gamma_n^{n+2}n^{n-2}}{2^{3n+1}(n+2)^{n+2}}\alpha_\mathcal{H}(\Omega)^{n(n+2)},
    \end{equation}
    where 
    \begin{equation*}
        \gamma_n:=\frac{\omega_{n-1}^n}{n^n\Big(\omega_n\big(4\pi^{n-2}\big)^{\frac1n}+2^{\frac{n-1}n}\omega_{n-1}\Big)^n}.
    \end{equation*}
\end{teorema}
We observe that \eqref{Tso}, if $f(x,u)=f(x)$, implies the following sharp estimate:
\begin{equation}\label{TsoABP}
    \|u\|_{L^\infty(\Omega)}\leq \frac{\mathrm{diam}(\Omega)}{2\omega_n^{1/n}}\left(\int_\Omega \det (-D^2u)dx\right)^{1/n},
\end{equation}
which is a version of the Alexandrov-Bakelman-Pucci (ABP) estimate with respect to the diameter.
In \cite{APABP}, a quantitative version of this estimate in terms of the area is proved in dimension two.
As an immediate consequence of Theorem \ref{Main Teo}, we obtain a quantitative version of \eqref{TsoABP}.
\begin{corollario}\label{QuantABP}
Under the assumptions of Theorem~\ref{Main Teo},
\begin{equation}\label{RisABP}
\frac{\dfrac{\mathrm{diam}(\Omega)}{2\,\omega_n^{1/n}}\,\|f_0\|_{L^1(\Omega)}^{\frac1n}-\|u\|_{L^\infty(\Omega)}}{\|u\|_{L^\infty(\Omega)}}
\ \geq\ \frac{\gamma_n^{n+2}n^{n-2}}{2^{3n+1}(n+2)^{n+2}}\,\alpha_\mathcal{H}(\Omega)^{n(n+2)} .
\end{equation}
\end{corollario}
We observe that Talenti-type comparison principles have been extended in several directions. On the one hand, different boundary conditions have been considered, including Neumann conditions \cite{CNT} and Robin conditions \cite{ANT} (see also \cite{AmatoBarbatoCitoMasielloPaoli} for a quantitative improvement). The comparison has also been extended to parabolic equations \cite{ALVINO199037}, nonlinear elliptic equations in divergence form \cite{BettaMercaldo, Talenti79}, fully nonlinear equations  not necessarily in divergence form \cite{SALANI2015763},
 the $p$-Laplace equation with Robin boundary conditions \cite{AGM}, $p$-Laplace equation in multiply connected domains \cite{BarbatoSalerno}, first-order Hamilton-Jacobi equations \cite{GN} (see \cite{AmatoBarbato1} for a quantitative improvement), anisotropic equations \cite{DPG,San} and nonlocal equations \cite{BrandoliniBFV, FeronePiscitelliVolzone,  FeroneVolzone21, FeroneVolzone23}. 

\medskip
The paper is organized as follows. In Section \ref{Sec2} we introduce
the geometric tools and Tso's symmetrization, and recall the properties of the Monge-Ampère operator needed in the sequel. In Section \ref{Sec3} we collect some useful lemmas that we will apply in Section \ref{Sec4}, where the main results are proved.

\section{Notation and Preliminaries}\label{Sec2}
In this section we provide the classical definitions and results that we need in the following. The reader can refer to \cite[Chapter 1]{S} for more details.
\begin{definizione}
    Let $\Omega, K\subset \R^n$ be two bounded convex sets. We define the \emph{Minkowski sum} $(+)$ as
    \begin{equation*}
        \Omega+K:=\{x+y \, : \, x\in \Omega, \, y\in K\}.
    \end{equation*}
\end{definizione} 

\begin{definizione}
    Let $E,F$ be two convex sets in $\mathbb{R}^n$. The \emph{Hausdorff distance} between $E$ and $F$ is defined as
\[
    d_{\mathcal{H}}(E,F) := \inf\left\{ \varepsilon>0 \, : \,  E \subset F + \varepsilon B , \, F \subset E + \varepsilon B  \right\},
\]
where $B$ is the unit ball centered at the origin and $F + \varepsilon B$ is the Minkowski sum.
\end{definizione}

\begin{definizione}
    \label{support}
Let $\Omega$ be a convex set. The \emph{support} function of $\Omega$ is defined as
\begin{equation*}
        h(\Omega, \theta):=\max_{x\in \overline{\Omega}} (x\cdot \theta), \qquad \theta\in \mathbb{S}^{n-1}.
    \end{equation*}
\end{definizione}
The support function allows us to define, for a convex set $\Omega$, the following geometrical quantities.
\begin{definizione}
 The \emph{width function} $w(\Omega, \cdot)$ of $\Omega$
is defined as
$$w(\Omega, \theta):= h(\Omega, \theta) + h(\Omega, -\theta) \quad \text{for } \theta\in \mathbb{S}^{n-1}.$$

The quantity $w(\Omega, \theta)$ is the thickness of $\Omega$ in the direction $\theta$ and it represents the distance between
the two support hyperplanes of $\Omega$ orthogonal to $\theta$. The maximum of the width function,
$$\mathrm{diam}(\Omega):=\max_{\theta\in \mathbb{S}^{n-1}} w(\Omega, \theta)$$
is the diameter of $\Omega$.

The mean value of the width function is called the \emph{mean width} and it is denoted by 
\begin{equation}\label{mean width}
    w(\Omega):= \frac{2}{n\omega_n} \int_{\mathbb{S}^{n-1}}h(\Omega, \theta)\, d\theta.
\end{equation}

The \emph{Steiner point} $s(\Omega)$ is defined via the vector-valued integral

\begin{equation}\label{SteinerP}
s(\Omega):=\frac{1}{\omega_n}\int_{\mathbb{S}^{n-1}}h(\Omega, \theta)\theta\, d\theta.
\end{equation}
\end{definizione}

\begin{definizione}
    The \emph{Steiner ball} $B_\Omega$ of a convex set $\Omega$ is the ball centered at the Steiner point of $\Omega$ with diameter equal to $w(\Omega)$.
\end{definizione}

\begin{oss}\label{Rmk hauss dist}
    By Definition \ref{support}, it is possible to characterize the Hausdorff distance as follows

    $$d_\mathcal{H}(E, F)= \norma{h(E,\cdot)-h(F,\cdot)}_{L^\infty}.$$
    Now it is clear that
    \begin{equation*}
        \abs{s(E)-s(F)}\le n d_\mathcal{H}(E,F).
    \end{equation*}
\end{oss} 
\begin{oss}
    \label{width/diam}
    Let us observe that it is possible to compare the mean width and the diameter of $\Omega$, for any convex set $\Omega\subset\mathbb{R}^n$, as the following bounds hold

    $$\frac{2\omega_{n-1}}{n\omega_n}\le\frac{w(\Omega)}{\mathrm{diam}(\Omega)}\le 1.$$

    The upper bound follows from the definitions of diameter and mean width, while the lower bound is proved in \cite{alexandrov}.
\end{oss}
If $\Omega\subset\R^n$ is a compact convex set with non-empty interior, we define its \emph{polar body} as 
\begin{equation*}
    \Omega^o=\{x\in\mathbb{R}^n:\langle x,y\rangle\leq 1\,\forall y\in \Omega\}.
\end{equation*}
If $x\in\Omega$, we denote by $\Omega^o_{x}$ the polar body $(\Omega-x)^o$. In the following, if $\Omega$ is an open bounded convex set, we write $\Omega^o$ for the polar of its closure.
\subsection{Quermassintegrals}
 For the content of this section, we will refer to \cite{S}. Let $\Omega \subset \R^n$ be a non-empty, bounded, convex set, let $B$ be the unit ball centered at the origin and $\rho > 0$. We can write the Steiner formula for the Minkowski sum $\Omega+ \rho B$ as
\begin{equation}
    \label{Steiner_formula}
    \abs{\Omega + \rho B} = \sum_{i=0}^n \binom{n}{i} W_i(\Omega) \rho^{i} .
\end{equation}
The coefficients $W_i(\Omega)$ are known in the literature as quermassintegrals of $\Omega$. In particular, $W_0(\Omega) = \abs{\Omega}$,  $nW_1(\Omega) = P(\Omega)$ and $W_n(\Omega) = \omega_n$ where $\omega_n$ is the measure of $B$.

Formula \eqref{Steiner_formula} can be generalized to every quermassintegral, obtaining
\begin{equation*}
    W_j(\Omega+\rho B) = \sum_{i=0}^{n-j} \binom{n-j}{i} W_{j+i}(\Omega) \rho^i, \qquad j=0, \ldots, n-1.
\end{equation*}

If $\Omega$ has $C^2$ boundary, the quermassintegrals can be written in terms of the principal curvatures of $\Omega$. More precisely, we can define the $j$-th mean curvature $H_j$ as the normalized elementary symmetric function of the principal curvatures $\kappa_1, \ldots, \kappa_{n-1}$ of $\partial \Omega$, i.e.
 \[
 H_0 = 1, \qquad \qquad H_j = \binom{n-1}{j}^{-1} \sum_{1 \leq i_1 < \ldots < i_j \leq n-1} \kappa_{i_1} \ldots \kappa_{i_j}, \qquad j = 1,\ldots,n-1.
 \]
 We observe that $H_1$ is the normalized \emph{mean curvature}, while $H_{n-1}$ is the normalized \emph{Gauss curvature}.
 Then, the quermassintegrals can be written as
 \begin{equation*}
     W_j(\Omega) = \frac{1}{n} \int_{\partial \Omega} H_{j-1} \, d \mathcal{H}^{n-1}, \qquad j = 1,\ldots, n-1.
 \end{equation*}

Furthermore, the Alexandrov-Fenchel inequalities hold (see, for instance, \cite{HCG})
\begin{equation}
    \label{Alexandrov_Fenchel_inequalities}
    \biggl( \frac{W_j(\Omega)}{\omega_n} \biggr)^{\frac{1}{n-j}} \geq \biggl( \frac{W_i(\Omega)}{\omega_n} \biggr)^{\frac{1}{n-i}}, \qquad 0 \leq i < j \leq n-1,
\end{equation}
where equality holds if and only if $\Omega$ is a ball. When $i=0$ and $j=1$, formula \eqref{Alexandrov_Fenchel_inequalities} reduces to the classical isoperimetric inequality, i.e.
\[
    P(\Omega) \geq n \omega_n^{\frac{1}{n}} \abs{\Omega}^{\frac{n-1}{n}}.
\]

\begin{definizione}
    Let $\Omega$ be a compact convex set. We define the \textit{k-th mean radius}, $k=0,\dots,n-1$, of $\Omega$ as
\begin{equation*}            \zeta_k(\Omega)=\left(\frac{W_k(\Omega)}{\omega_n}\right)^\frac{1}{n-k}.
\end{equation*}
\end{definizione}
Then the Alexandrov-Fenchel inequalities can be rewritten in terms of the mean radii, obtaining
\begin{equation*}
    \zeta_i(\Omega)\leq\zeta_j(\Omega),\quad 0\leq i\leq j\leq n-1.
\end{equation*}
\begin{oss}
    We observe that, since $w(\Omega)=2\zeta_{n-1}(\Omega)$, Remark \ref{width/diam} implies 
    
    \begin{equation}
        \label{zeta/diam}
        \frac{\zeta_{n-1}(\Omega)}{\mathrm{diam}(\Omega)}\geq \frac{\omega_{n-1}}{n\omega_n}.
    \end{equation}
\end{oss}
\begin{oss}
    If $K$ is a $j$-dimensional convex set in $\R^n$, then the $i$-th quermassintegral, $i=0,\dots,n-j-1$, equals zero, while the others may be non-null.
\end{oss}

\begin{definizione}\label{Hauss}
The \textit{Hausdorff asymmetry index} is defined as 
\begin{equation*}
    \alpha_\mathcal{H}(\Omega)=\dfrac{d_\mathcal{H}(\Omega,B_\Omega)}{\zeta_{n-1}(\Omega)},
\end{equation*}
 where $B_\Omega$ is the Steiner ball of $\Omega$.  
\end{definizione}
We observe that, since the support function of $\Omega$ centered at $s(\Omega)$ is bounded by $\mathrm{diam}(\Omega)$, we have that
\begin{equation}\label{BoundHaus}
    \alpha_\mathcal{H}(\Omega)\leq\frac{\mathrm{diam}(\Omega)}{\zeta_{n-1}(\Omega)}\leq\frac{n\omega_n}{\omega_{n-1}},
\end{equation}
where the last inequality follows by \eqref{zeta/diam}.
\subsection{Symmetrization}
Let $\Omega\subset\mathbb R^n$ be a measurable set with
$0<|\Omega|<\infty$ and let $u:\Omega\to\mathbb R$ be a
measurable function. For $t\ge0$, set
\[
    \mu(t):=
    \bigl|\{x\in\Omega:|u(x)|>t\}\bigr|.
\]
The one-dimensional \emph{decreasing rearrangement} of $u$ is defined by
\[
    u^*(s):=\inf\{t\ge0:\mu(t)\le s\},
    \qquad 0< s\le|\Omega|,\quad u^*(0)=\|u\|_{L^\infty(\Omega)}.
\]
The \emph{Schwarz rearrangement} of $u$ is defined by
\[
    u^\#_0(x):=u^*(\omega_n|x|^n),
    \qquad x\in\Omega^*_0,
\]
where $\Omega^*_0$ is the ball centered at the origin with
$|\Omega^*_0|=|\Omega|$. 
The functions $|u|$, $u^*$ and $u^\#_0$ are equimeasurable. Moreover, the Hardy-Littlewood inequality gives
\begin{equation*}
    \int_E |u(x)|\,dx
    \le
    \int_0^{|E|}u^*(s)\,ds
\end{equation*}
for every measurable set $E\subseteq\Omega$.\\
Now let $\Omega$ be an open, bounded and convex set in $\mathbb{R}^n$ and let
\[
\mathcal A(\Omega)=\big\{u\in C^2(\Omega)\cap C(\overline\Omega):\ u \text{ is strictly concave in }\Omega
\text{ and vanishes on }\partial\Omega\big\}.
\]
We observe that any $u\in\mathcal A(\Omega)$ is positive in $\Omega$ and attains its
maximum $M:=\|u\|_{L^\infty(\Omega)}$ at a unique point of $\Omega$. Denote
    \begin{equation*}
        \Omega_t=\{x\in \Omega\,:\, u(x)>t\},\quad\partial\Omega_t=\{x\in\Omega\,:\,u(x)=t\},
    \end{equation*}
    for $t\in[0,M]$.
In \cite{tso1989}, K. Tso extended the Schwarz rearrangement, defining the $k$-symmetrization of $u$, $k=0,\dots,n-1$, as
\[
u^\#_k(x)=\sup\big\{t\in[0,M]:\ \zeta_k(\Omega_t)>|x|\big\},
\]
which is a radially symmetric function defined on the ball
$\Omega^*_k=\{x:|x|<\zeta_k(\Omega)\}$ (note that $W_k(\Omega^*_k)=W_k(\Omega)$).
\\ 
It is possible to prove that $u^\#_k$ is
strictly decreasing along the radii and strictly concave, so that
$u^\#_k\in\mathcal A(\Omega^*_k)$.

\subsection{The Monge-Ampère operator}
Let $\Omega$ be an open subset of $\R^n$ and let $u\in C^2(\Omega)$. The Monge-Ampère operator is fully nonlinear and it is non-elliptic, unless one restricts to the class  of concave functions (see \cite{DPF, gutierrez2016} for more details on the topic).
We recall that the operator $\det^{\frac{1}{n}}$ is homogeneous of degree 1, and if we denote by

$$S_n^{ij}(-D^2u)=\frac{\partial}{\partial u_{ij}}\det(-D^2u),$$
the Euler identity for homogeneous functions gives
\begin{equation*}
    \det(-D^2u)=\frac{1}{n} S_n^{ij}(-D^2u)u_{ij}.
\end{equation*}
A direct computation shows that the $\left(S_n^{1j}(-D^2u), \dots, S_n^{nj}(-D^2u)\right)$ is divergence-free, hence $\det(-D^2u)$ can be written in divergence form

\begin{equation}
    \label{skdiv}
    \det(-D^2u)=\frac{1}{n}\left(S_n^{ij}(-D^2u) u_j\right)_i,
\end{equation}
where the subscripts $i, j$ stand for partial differentiation. 

If $u\in C^2(\Omega)$ and $t$ is a regular value of $u$, on the boundary of $\partial\Omega_t$ it is possible to link the Monge-Ampère operator with the Gauss curvature

\begin{equation}
    \label{hksk}
    H_{n-1}=-\frac{S_n^{ij}(-D^2u)u_iu_j}{\abs{\nabla u}^{n+1}}.
\end{equation}
We recall the following classical result without proof (see \cite[Proposition 6]{tso1989} for details).
\begin{prop}
    \label{prop6}
    Let $\Sigma(t)$ be a smooth family of embedded hypersurfaces in $\mathbb{R}^n$ oriented with respect to the unit inner normal $-\nu$. Then 
    \begin{equation*}
        \frac{d}{dt}\int_{\Sigma(t)}H_{n-2}\, d\mathcal{H}^{n-1}=\int_{\Sigma(t)}H_{n-1}\,\nu\cdot\xi \, d\mathcal{H}^{n-1}
    \end{equation*}
    where 
    \begin{equation*}
        \xi=\frac{\partial\Sigma(t)}{\partial t}
    \end{equation*}
    is the variation vector field.
\end{prop}
We recall that, if $u\in C^2(\Omega)$ is a solution to \eqref{MAequation}, then it is strictly concave. Moreover, for $0<t<\|u\|_{L^\infty(\Omega)}$, writing $\zeta_{n-1}(t)=\zeta_{n-1}(\Omega_t)$,
we have the following identities
\begin{equation}\label{HQ:identities}
\int_{\Omega_t}\det(-D^2u)\,dx=\frac{1}{n}\int_{\partial\Omega_t}H_{n-1}|\nabla u|^n\,d\mathcal H^{n-1},\qquad
\int_{\partial\Omega_t}\frac{H_{n-1}}{|\nabla u|}\,d\mathcal H^{n-1}=n\omega_n\,(-\zeta'_{n-1}(t)),
\end{equation}
where the first identity follows from \eqref{skdiv}-\eqref{hksk} and the divergence Theorem, while the second follows from Proposition \ref{prop6}. Moreover, applying H\"older's inequality to the second identity in \eqref{HQ:identities} and using the first identity in \eqref{HQ:identities}, we have
\begin{equation}\label{TalentiIneq}
\omega_n\le(-\zeta'_{n-1}(t))^n\int_{\Omega_t}\det(-D^2u)\,dx .
\end{equation}
We observe that, from the concavity of $u$, for $t,\,s\in[0,M)$ and $\lambda\in[0,1]$, we have $(1-\lambda)\Omega_t+\lambda\Omega_s\subset\Omega_{(1-\lambda)t+\lambda s}$. Then from the Brunn-Minkowski inequality for $W_{n-1}$ (see, for instance, \cite{S}) and the monotonicity of the quermassintegrals with respect to inclusion, we have that $\zeta_{n-1}(t)$ is a concave and decreasing function on $[0,\|u\|_{L^\infty(\Omega)})$. From the concavity we have
\begin{equation*}
    \zeta_{n-1}(\tau)\leq\zeta_{n-1}(t)+\zeta'_{n-1}(t)(\tau-t),
\end{equation*}
and, as $\tau\rightarrow\|u\|_{L^\infty(\Omega)}$, since $\zeta_{n-1}(\tau)\rightarrow0$, we have
\begin{equation}\label{HQ:concave}
\frac{\zeta_{n-1}(t)}{-\zeta'_{n-1}(t)}\ \ge\ \|u\|_{L^\infty(\Omega)}-t
\qquad\text{for every }0\le t<\|u\|_{L^\infty(\Omega)} .
\end{equation}
We also recall that the solution to \eqref{MAequationSymm} is given by (see \cite{tso1989} for the proof)
\begin{equation}\label{v radiale}
    v(r)=\int_r^{\zeta_{n-1}(\Omega)}\left(\frac1{\omega_n}\int_0^{\omega_n\rho^n}f^*(s)\,ds\right)^{\frac1n}d\rho,
\qquad 0\le r\le\zeta_{n-1}(\Omega).
\end{equation}
\begin{oss}
The function $v-u^\#_{n-1}$ is radially nonincreasing, so that
\begin{equation}\label{normaTso}
\|v-u^\#_{n-1}\|_{L^\infty(\Omega^*_{n-1})}=v(0)-\|u\|_{L^\infty(\Omega)} .
\end{equation}
Indeed, let $\varphi(r):=v(r)-u^\#_{n-1}(r)$ be its radial profile. 
By the explicit expression \eqref{v radiale} of $v$, we have
\begin{equation*}
    v'(r)=-\Big(\frac{1}{\omega_n}\int_0^{\omega_nr^n} f^*(s)\,ds\Big)^{\frac1n},
\qquad 0\le r<\zeta_{n-1}(\Omega).
\end{equation*}
Let now $t\in(0,\|u\|_{L^\infty(\Omega)})$ and set $r=\zeta_{n-1}(t)$. By definition of $u^\#_{n-1}$, we have 
$u^\#_{n-1}(r)=t$, hence
\begin{equation*}
\frac{d}{dr}u^\#_{n-1}(r)=\frac1{\zeta'_{n-1}(t)}=-\frac1{-\zeta'_{n-1}(t)} .
\end{equation*}
Therefore
\begin{equation*}
    \varphi'(r)=-\Big(\frac{1}{\omega_n}\int_0^{\omega_n\zeta_{n-1}(t)^n} f^*(s)\,ds\Big)^{\frac1n}+\frac1{-\zeta'_{n-1}(t)} .
\end{equation*}
By \eqref{TalentiIneq}, the monotonicity of $f$, the Hardy-Littlewood inequality and 
$|\Omega_t|\le\omega_n\zeta_{n-1}(t)^n$, we have
\begin{equation*}
\frac{\omega_n}{(-\zeta'_{n-1}(t))^n}\le\int_{\Omega_t}\det(-D^2u)\,dx\le\int_{\Omega_t}f_0\,dx
\le\int_0^{\omega_n\zeta_{n-1}(t)^n}f^*(s)\,ds ,
\end{equation*}
which is exactly Tso's comparison at the level of derivatives, that is, $\varphi'(r)\le0$. 
Hence $\varphi$ is nonincreasing and attains its maximum at $r=0$; since 
$u^\#_{n-1}(0)=\|u\|_{L^\infty(\Omega)}$, \eqref{normaTso} follows.
\end{oss}
\subsection{The support function of the level sets}\label{SupFunLevelSet}
Let $u$ be a solution to \eqref{MAequation} and let $x_0$ be its maximum point.
For $\theta\in\mathbb{S}^{n-1}$ and $0\le t<\|u\|_{L^\infty(\Omega)}$, let
\[
h(\theta,t):=\sup_{y\in\Omega_t}\langle y-x_0,\theta\rangle
\]
be the support function of $\overline{\Omega_t}$ with respect to $x_0$, and write
$h_\Omega(\theta):=h(\theta,0)$. The map
\begin{equation*}
    x\mapsto\left(-\frac{\nabla u(x)}{\abs{\nabla u(x)}},\,u(x)\right)
\end{equation*}
is a $C^1$ diffeomorphism of
$\Omega\setminus\{x_0\}$ onto $\mathbb{S}^{n-1}\times(0,\|u\|_{L^\infty(\Omega)})$; we
denote by $x(\theta,t)$ its inverse, that is, the unique point of $\partial\Omega_t$ where
the outer unit normal to $\Omega_t$ is $\theta$. Then
\[
h(\theta,t)=\langle x(\theta,t)-x_0,\theta\rangle,\qquad
\nabla u(x(\theta,t))=-|\nabla u(x(\theta,t))|\,\theta .
\]
Differentiating $u(x(\theta,t))=t$ with respect to $t$ gives
$\langle\nabla u(x(\theta,t)),\partial_tx(\theta,t)\rangle=1$, hence
\[
\partial_th(\theta,t)=\langle\partial_tx(\theta,t),\theta\rangle
=-\frac1{|\nabla u(x(\theta,t))|}.
\]
Integrating on $[0,T]$ and letting $T\to\|u\|_{L^\infty(\Omega)}$ (recall that
$|h(\theta,T)|\le \mathrm{diam}(\Omega_T)\to0$), we obtain
\begin{equation}\label{h e gradiente}
h_\Omega(\theta)=\int_0^{\|u\|_{L^\infty(\Omega)}}\frac{d\tau}{|\nabla u(x(\theta,\tau))|}.
\end{equation}
Since 
\begin{equation*}
    \zeta_{n-1}(\Omega)=\int_0^{\|u\|_{L^\infty(\Omega)}}(-\zeta'_{n-1}(\tau))\,d\tau,
\end{equation*}
from \eqref{h e gradiente} we have
\[
h_\Omega(\theta)-\zeta_{n-1}(\Omega)
=\int_0^{\|u\|_{L^\infty(\Omega)}}\left(\frac1{|\nabla u(x(\theta,\tau))|}+\zeta'_{n-1}(\tau)\right)d\tau .
\]
Then by Fubini-Tonelli and the change of variables given by the Gauss map of
$\partial\Omega_\tau$, whose Jacobian is $H_{n-1}$, we find
\begin{equation}\label{HQ:supp1}
\int_{\mathbb{S}^{n-1}}\big|h_\Omega-\zeta_{n-1}(\Omega)\big|\,d\mathcal H^{n-1}
\le\int_0^{\|u\|_{L^\infty(\Omega)}}\int_{\partial\Omega_\tau}H_{n-1}
\left|\frac1{|\nabla u|}+\zeta'_{n-1}(\tau)\right|d\mathcal H^{n-1}\,d\tau .
\end{equation}
We finally recall that \eqref{mean width} gives
\begin{equation}\label{ZetaSupport}
    \frac1{n\omega_n}\int_{\mathbb{S}^{n-1}}h_\Omega\,d\mathcal H^{n-1}=\zeta_{n-1}(\Omega).
\end{equation}
\section{Useful Lemmas}\label{Sec3}
In this section we list some lemmas that we will apply below. The first result is a quantitative improvement of H\"older's inequality. This result has been proved in \cite{Aldaz} and the version given below is only a small modification of it, rewritten in a form that will be more useful in the following.
\begin{lemma}
    Let $2\leq p<+\infty$, $p'=\frac{p}{p-1}$ and let $f\in L^p$ and $g\in L^{p'}$ satisfy $\|f\|_{L^p},\,\|g\|_{L^{p'}}>0$. Then
    \begin{equation}\label{AldazL1}
        \|fg\|_{L^1}\leq \|f\|_{L^p}\|g\|_{L^{p'}}\left(1-\frac{1}{4p}\left\|\frac{|f|^p}{\|f\|_{L^p}^p}-\frac{|g|^{p'}}{\|g\|_{L^{p'}}^{p'}} \right\|_{L^1}^2 \right).
    \end{equation}
\end{lemma}
\begin{proof}
We recall that the inequality proved in \cite[Theorem 2.2]{Aldaz} is
\begin{equation}\label{Aldaz}
        \|fg\|_{L^1}\leq \|f\|_{L^p}\|g\|_{L^{p'}}\left(1-\frac{1}{p}\left\|\frac{|f|^\frac{p}{2}}{\|f\|_{L^p}^\frac{p}{2}}-\frac{|g|^\frac{p'}{2}}{\|g\|_{L^{p'}}^\frac{p'}{2}} \right\|_{L^2}^2 \right).
    \end{equation}
    Set
    \begin{equation*}
        F=\frac{|f|^p}{\|f\|_{L^p}^p} \quad\text{ and }\quad G=\frac{|g|^{p'}}{\|g\|_{L^{p'}}^{p'}}.
    \end{equation*}
    Minkowski's inequality gives
    \begin{equation*}
        \|\sqrt{F}+\sqrt{G}\|_{L^2}\leq\|F\|_{L^1}^\frac{1}{2}+\|G\|_{L^1}^\frac{1}{2}=2.
    \end{equation*}
    Then, combining the previous estimate with H\"older's inequality, we get
    \begin{equation}\label{L^2 to L^1}
        \|F-G\|_{L^1}=\left\|\left(\sqrt{F}-\sqrt{G}\right)\left(\sqrt{F}+\sqrt{G}\right)\right\|_{L_1}\leq\left\|\sqrt{F}-\sqrt{G}\right\|_{L^2}\left\|\sqrt{F}+\sqrt{G}\right\|_{L^2}\leq 2\left\|\sqrt{F}-\sqrt{G}\right\|_{L^2}.
    \end{equation}
    Finally, combining \eqref{Aldaz}-\eqref{L^2 to L^1}, we have
    \begin{equation*}
        \|fg\|_{L^1}\leq \|f\|_{L^p}\|g\|_{L^{p'}}\left(1-\frac{1}{4p}\left\|\frac{|f|^p}{\|f\|_{L^p}^p}-\frac{|g|^{p'}}{\|g\|_{L^{p'}}^{p'}} \right\|_{L^1}^2 \right).
    \end{equation*}
\end{proof}

The next lemma is a generalization of the planar estimate proved by Talenti in \cite{talenti1981}. 
\begin{lemma}
    Let $\Omega\subset\R^n$ be an open bounded convex set and let $u\in C^2(\Omega)\cap C(\overline{\Omega})$ be a solution to \eqref{MAequation} attaining its maximum at $x_0\in\Omega$. Then   \begin{equation}\label{boundNormaInf di u}
        \|u\|_{L^\infty(\Omega)}\leq\frac{2n}{\omega_n^\frac{1}{n}}\zeta_0(\Omega)\|f_0\|_{L^1(\Omega)}^\frac{1}{n}.
    \end{equation}
\end{lemma}
\begin{proof}
    We recall the following a priori estimate (see \cite{APABP})
    \begin{equation*}
        \|u\|_{L^\infty(\Omega)}^n\leq \frac{1}{\abs{\Omega^o_{x_0}}}\int_\Omega\det (-D^2u)\,dx.
    \end{equation*}
    Using that $\det(-D^2u)=f(x,u)\leq f_0(x)$ in the latter, we get
    \begin{equation}\label{A priori AP}
        \|u\|_{L^\infty(\Omega)}^n\leq\frac{\abs{\Omega}}{\abs{\Omega^o_{x_0}}\abs{\Omega}}\|f_0\|_{L^1(\Omega)}.
    \end{equation}
    By John's Theorem, there exists a unique ellipsoid $E$ of maximal volume contained in $\Omega$, with center $c$, such that $E\subset \Omega\subset c+n(E-c)$.
    Since $x_0\in\Omega\subset c+n(E-c)$, then $x_0-c\in n(E-c)$ and then $c-x_0\in n(E-c)$, since $E-c$ is symmetric about the origin. Then
    \begin{equation*}
        \Omega-x_0\subset (c-x_0)+n(E-c)\subset 2n(E-c).
    \end{equation*}
    Passing to the polar body, the latter gives
    \begin{equation*}
        \Omega^o_{x_0}=(\Omega-x_0)^o \supset\frac{1}{2n}(E-c)^o.
    \end{equation*}
    Then
    \begin{equation}\label{stima Polare}
        \abs{\Omega^o_{x_0}}\geq\frac{1}{(2n)^n}\abs{(E-c)^o}=\frac{\omega_n^2}{(2n)^n\abs{E}}\geq\frac{\omega_n^2}{(2n)^n\abs{\Omega}},
    \end{equation}
    where in the last inequality we have used $E\subset\Omega$. Combining \eqref{A priori AP}-\eqref{stima Polare} we get \eqref{boundNormaInf di u}.
\end{proof}
We observe that, if we consider $n=2$ in \eqref{boundNormaInf di u}, we do not recover the same constant as in \cite{talenti1981}. The reason is that Talenti applied the planar Mahler inequality in \eqref{A priori AP}, however this inequality is no longer available for general dimension, then we use a rougher estimate provided by John's Theorem.\\
We state the following elementary lemma without proof.
\begin{lemma}\label{f decrescente}
    Let $F:[0,+\infty)\rightarrow\R$ be a nonincreasing and locally integrable function and consider real numbers $0<a<b$. Then
    \begin{equation*}
        \frac{1}{a}\int_0^a F(x)\,dx\geq\frac{1}{b}\int_0^b F(x)\,dx.
    \end{equation*}
\end{lemma}
Finally we prove the following lemma, which relates the difference between the support function of $\Omega$ and $\zeta_{n-1}(\Omega)$ with the Hausdorff asymmetry.
\begin{lemma}\label{lem:geom}
Let $\Omega\subset\mathbb{R}^n$ be an open, bounded and convex set, let $x_0\in\Omega$,
and let $h_\Omega=h(\Omega-x_0,\cdot)$ be the support function of $\Omega$ with respect
to $x_0$. Then
\begin{equation}\label{Supp0}
\frac1{n\omega_n}\int_{\mathbb{S}^{n-1}}\big|h_\Omega-\zeta_{n-1}(\Omega)\big|\,d\mathcal H^{n-1}
\ \ge\ \gamma_n\,\zeta_{n-1}(\Omega)\,\alpha_{\mathcal H}(\Omega)^n ,
\end{equation}
where $\gamma_n$ is the constant in Theorem~\ref{Main Teo}.
\end{lemma}
\begin{proof}
Let
\[
m:=\max_{\theta\in\mathbb{S}^{n-1}}\big|h_\Omega(\theta)-\zeta_{n-1}(\Omega)\big|
=d_{\mathcal H}\big(\Omega,{B_{\zeta_{n-1}(\Omega)}(x_0)}\big),
\]
where the second equality follows by Remark \ref{Rmk hauss dist}, and let $\bar\theta\in\mathbb{S}^{n-1}$ be a point where the maximum is attained.
Since $x_0\in\Omega$, we have $0\le h_\Omega\le \mathrm{diam}(\Omega)$. Moreover
$2\zeta_{n-1}(\Omega)\le \mathrm{diam}(\Omega)$, so that 
\begin{equation*}
    -\frac{\mathrm{diam}(\Omega)}{2}\leq h_\Omega(\theta)-\zeta_{n-1}(\Omega)\leq\mathrm{diam}(\Omega),
\end{equation*}
that is $\abs{h_\Omega(\theta)-\zeta_{n-1}(\Omega)}\leq \mathrm{diam}(\Omega)$ for all $\theta\in\mathbb{S}^{n-1}$. Then $m\leq\mathrm{diam}(\Omega)$.
Let
\[
\rho:=\frac{m}{2\mathrm{diam}(\Omega)}\le\frac12,
\qquad
\Sigma:=\left\{\theta\in\mathbb{S}^{n-1}:\ d_{\mathbb{S}^{n-1}}(\theta,\bar\theta)\leq\rho\right\},
\]
where $d_{\mathbb{S}^{n-1}}$ is the geodesic distance. We recall that $h_\Omega$ is Lipschitz, with $\mathrm{Lip}(h_\Omega)\leq\mathrm{diam}(\Omega)$, that is,
\begin{equation*}
    \abs{h_\Omega(\theta)-h_\Omega(\theta')}\leq\mathrm{diam}(\Omega)\abs{\theta-\theta'}, \quad\text{for all }\theta,\,\theta'\in\mathbb{S}^{n-1}.
\end{equation*}
Then, if $\theta\in\Sigma$, we have $\abs{\theta-\overline{\theta}}\leq\rho$ and then $\abs{h_\Omega(\theta)-h_\Omega(\overline{\theta})}\leq m/2$. This gives
\begin{equation}\label{per dopo}
    \abs{h_\Omega(\theta)-\zeta_{n-1}(\Omega)}\geq\abs{h_\Omega(\overline{\theta})-\zeta_{n-1}(\Omega)}-\abs{h_\Omega(\theta)-h_\Omega(\overline{\theta})}\geq \frac{m}{2},
\end{equation}
where we have used that the maximum $m$ is attained at $\overline{\theta}$.
Using spherical coordinates centered at $\bar\theta$, in which
$d\mathcal H^{n-1}=(\sin s)^{n-2}\,ds\,d\mathcal H^{n-2}$, where $s$ is the geodesic distance from $\bar\theta$,
and $\mathcal H^{n-2}(\mathbb S^{n-2})=(n-1)\omega_{n-1}$, we have
\[
\mathcal H^{n-1}(\Sigma)=(n-1)\,\omega_{n-1}\int_0^\rho\sin^{n-2}s\,ds.
\]
Since $\rho\leq1/2$ and $\sin s\geq(2/\pi)s$ on $[0,\pi/2]$, the latter gives
\begin{equation*}
    \mathcal H^{n-1}(\Sigma)\geq\omega_{n-1}\left(\frac2\pi\right)^{n-2}\rho^{n-1}=\frac{\omega_{n-1}m^{n-1}}{2\pi^{n-2}\mathrm{diam}(\Omega)^{n-1}}.
\end{equation*}
Then from \eqref{per dopo}
\[
\int_{\mathbb{S}^{n-1}}\big|h_\Omega-\zeta_{n-1}(\Omega)\big|\,d\mathcal H^{n-1}
\ge\frac m2\,\mathcal H^{n-1}(\Sigma)
\ge\frac{\omega_{n-1}m^n}{4\pi^{n-2}\mathrm{diam}(\Omega)^{n-1}}.
\]
From the definition of $m$, the latter is equivalent to
\begin{equation}\label{HQ:Linf}
d_{\mathcal H}\big(\Omega,{B_{\zeta_{n-1}(\Omega)}(x_0)}\big)
\le\left(\frac{4\pi^{n-2}}{\omega_{n-1}}\,\mathrm{diam}(\Omega)^{n-1}
\int_{\mathbb{S}^{n-1}}\big|h_\Omega-\zeta_{n-1}(\Omega)\big|\,d\mathcal H^{n-1}\right)^{\frac1n}.
\end{equation}
Let $s(\Omega)$ be the Steiner point of $\Omega$ defined as in \eqref{SteinerP}.
Since $\int_{\mathbb{S}^{n-1}}\theta\,d\mathcal H^{n-1}=0$, we get
\[
|s(\Omega)-x_0|=\frac1{\omega_n}\left|\int_{\mathbb{S}^{n-1}}\big(h_\Omega(\theta)-\zeta_{n-1}(\Omega)\big)\theta\,d\mathcal H^{n-1}\right|
\le\frac1{\omega_n}\int_{\mathbb{S}^{n-1}}\big|h_\Omega-\zeta_{n-1}(\Omega)\big|\,d\mathcal H^{n-1}.
\]
Then, from the triangle inequality and \eqref{ZetaSupport},
we get $|s(\Omega)-x_0|\le2n\,\zeta_{n-1}(\Omega)$. Writing
\begin{equation*}
    |s(\Omega)-x_0|=|s(\Omega)-x_0|^{\frac1n}|s(\Omega)-x_0|^{\frac{n-1}n},
\end{equation*}
 we obtain
\begin{equation}\label{HQ:steiner}
|s(\Omega)-x_0|\le\left(\frac1{\omega_n}\int_{\mathbb{S}^{n-1}}\big|h_\Omega-\zeta_{n-1}(\Omega)\big|\,d\mathcal H^{n-1}\right)^{\frac1n}
\big(2n\,\zeta_{n-1}(\Omega)\big)^{\frac{n-1}n}.
\end{equation}
Since $B_\Omega=B_{\zeta_{n-1}(\Omega)}(s(\Omega))$, we have
\[
d_{\mathcal H}(\Omega,B_\Omega)\le
d_{\mathcal H}\big(\Omega,B_{\zeta_{n-1}(\Omega)}(x_0)\big)+|s(\Omega)-x_0| .
\]
Combining \eqref{HQ:Linf}, \eqref{zeta/diam} and
\eqref{HQ:steiner}, we get
\[
\alpha_\mathcal{H}(\Omega)=\frac{d_{\mathcal H}(\Omega,B_\Omega)}{\zeta_{n-1}(\Omega)}\le
n\left(\frac{\omega_n}{\omega_{n-1}}\big(4\pi^{n-2}\big)^{\frac1n}+2^{\frac{n-1}n}\right)
\left(\frac1{n\omega_n\,\zeta_{n-1}(\Omega)}\int_{\mathbb{S}^{n-1}}\big|h_\Omega-\zeta_{n-1}(\Omega)\big|\,d\mathcal H^{n-1}\right)^{\frac1n}.
\]
Raising to the power $n$, we get \eqref{Supp0}.
\end{proof}
\section{Proofs of the main results}\label{Sec4}
The proof of Theorem \ref{Main Teo} requires several steps, which we organize into two lemmas. Below we give a sketch of the proof, in order to make the section more readable. \\
Let us suppose that the assumptions of Theorem \ref{Main Teo} hold true.
The first step is an application of the quantitative H\"older inequality proved in \cite{Aldaz}, that we will apply instead of the classical H\"older inequality used in \cite{talenti1981, tso1989}. This will give us a quantitative improvement with respect to the following deficit
\begin{equation}\label{HQ:defD}
\mathcal D(t):=\bigint_{\partial\Omega_t}\left|
\frac{H_{n-1}|\nabla u|^n}{\displaystyle\int_{\partial\Omega_t}H_{n-1}|\nabla u|^n\,d\mathcal H^{n-1}}
-\frac{\dfrac{H_{n-1}}{|\nabla u|}}{\displaystyle\int_{\partial\Omega_t}\frac{H_{n-1}}{|\nabla u|}\,d\mathcal H^{n-1}}
\right|d\mathcal H^{n-1},\quad0<t<\|u\|_{L^\infty(\Omega)}.
\end{equation}
This deficit measure \emph{how far} $\abs{\nabla u}$ is from being \emph{constant} on the level surface $\partial\Omega_t$. This will be the content of Lemma \ref{HQ:lemma}. We then relate this \emph{H\"older deficit} to a suitable power of the Hausdorff asymmetry index. This is achieved by rewriting each quantity that appears in the integral in terms of the support function (as explained in Subsection \ref{SupFunLevelSet}) and applying Lemma \ref{lem:geom}. This will be the content of Lemma \ref{HQ:thm}. Finally we will be ready for the proof of Theorem \ref{Main Teo}.\\
We start with the following lemma.
\begin{lemma}\label{HQ:lemma}
Let $v$ be the solution to \eqref{MAequationSymm}. Then, under the assumptions of Theorem \ref{Main Teo}, for every $0\le\tau<\|u\|_{L^\infty(\Omega)}$,
\begin{equation}\label{HQ:31}
\frac{\omega_n^{\frac{n-1}{n}}}{4n\,\|f_0\|_{L^1(\Omega)}^{\frac{n-1}{n}}}
\int_0^\tau\frac{\mathcal D(t)^2}{(-\zeta'_{n-1}(t))^{n-1}}\,dt\ \le\ v(\zeta_{n-1}(\tau))-\tau,
\end{equation}
where $\mathcal D(t)$ is defined in \eqref{HQ:defD}.
In particular,
\begin{equation}\label{HQ:HD}
\|v-u^\#_{n-1}\|_{L^\infty(\Omega^*_{n-1})}\ \ge\
\frac{\omega_n^{\frac{n-1}{n}}}{4n\,\|f_0\|_{L^1(\Omega)}^{\frac{n-1}{n}}}
\int_0^{\|u\|_{L^\infty(\Omega)}}\frac{\mathcal D(t)^2}{(-\zeta'_{n-1}(t))^{n-1}}\,dt .
\end{equation}
\end{lemma}
\begin{proof}
We recall that $u\ge0$ in $\Omega$ and $f(x,u(x))\le f_0(x)$.

\smallskip
\noindent\emph{Step 1: quantitative estimate on the level sets.}
Let $0<t<\|u\|_{L^\infty(\Omega)}$. We apply the quantitative H\"older
inequality \eqref{AldazL1} on $\partial\Omega_t$, with $p=n+1$ and $p'=\frac{n+1}{n}$, to the functions
\[
f=H_{n-1}^{\frac1{n+1}}\,|\nabla u|^{\frac n{n+1}}
\qquad\text{and}\qquad
g=H_{n-1}^{\frac n{n+1}}\,|\nabla u|^{-\frac n{n+1}},
\]
obtaining
\begin{equation*}
    \int_{\partial\Omega_t}H_{n-1}\,d\mathcal{H}^{n-1}\leq\left(\int_{\partial\Omega_t}H_{n-1}\abs{\nabla u}^n\,d\mathcal{H}^{n-1}\right)^\frac{1}{n+1}\left(\int_{\partial\Omega_t}\frac{H_{n-1}}{\abs{\nabla u}}\,d\mathcal{H}^{n-1}\right)^\frac{n}{n+1}\left[1-\frac{\mathcal{D}(t)^2}{4(n+1)}\right].
\end{equation*}

We observe that, since each of the two quotients
appearing in \eqref{HQ:defD} has integral equal to $1$ on $\partial\Omega_t$, we have
$0\le\mathcal D(t)\le2$. Using \eqref{HQ:identities} and the Gauss-Bonnet Theorem, we can rewrite the H\"older inequality above as
\[
n\omega_n\le\left(n\int_{\Omega_t}\det(-D^2u)\,dx\right)^{\frac1{n+1}}
\big(n\omega_n(-\zeta'_{n-1}(t))\big)^{\frac n{n+1}}
\left(1-\frac{\mathcal D(t)^2}{4(n+1)}\right),
\]
that is,
\[
\omega_n\le(-\zeta'_{n-1}(t))^n\left(\int_{\Omega_t}\det(-D^2u)\,dx\right)
\left(1-\frac{\mathcal D(t)^2}{4(n+1)}\right)^{n+1}.
\]
Using $(1-s)^{-(n+1)}\ge1+(n+1)s$ for $0\le s<1$ with
$s=\mathcal D(t)^2/(4(n+1))\le1/(n+1)$, we obtain
\begin{equation}\label{quant1}
\int_{\Omega_t}\det(-D^2u)\,dx\ \ge\ \frac{\omega_n}{(-\zeta'_{n-1}(t))^n}\left(1+\frac{\mathcal D(t)^2}{4}\right).
\end{equation}
The latter can be seen as a quantitative
improvement of \eqref{TalentiIneq}.

\smallskip
\noindent\emph{Step 2: use of the equation.}
Set $\mu(t)=\abs{\Omega_t}$. The monotonicity of $f$, the Hardy-Littlewood inequality and
$\mu(t)\le\omega_n\zeta_{n-1}(t)^n$ (that is, $\zeta_0(\Omega_t)\le\zeta_{n-1}(\Omega_t)$)
give
\[
\int_{\Omega_t}\det(-D^2u)\,dx\le\int_{\Omega_t}f_0\,dx\le\int_0^{\mu(t)}f^*(s)\,ds\le
\int_0^{\omega_n\zeta_{n-1}(t)^n}f^*(s)\,ds .
\]
Set
\[
g_1:=\int_0^{\omega_n\zeta_{n-1}(t)^n}f^*(s)\,ds,\qquad
g_2:=\frac{\omega_n\,\mathcal D(t)^2}{4(-\zeta'_{n-1}(t))^n}.
\]
By \eqref{quant1}, $\omega_n(-\zeta'_{n-1}(t))^{-n}\le g_1-g_2$; in particular
$0\le g_2<g_1\le\|f_0\|_{L^1(\Omega)}$. Hence
\[
1\le\frac{-\zeta'_{n-1}(t)}{\omega_n^{1/n}}\,(g_1-g_2)^{\frac1n}.
\]
Combining the latter and the following elementary inequality, which follows from the concavity of
$s\mapsto s^{1/n}$ and from $g_1\le\|f_0\|_{L^1(\Omega)}$,
\begin{equation*}
(g_1-g_2)^{\frac1n}\le g_1^{\frac1n}-\frac{g_2}{n\,\|f_0\|_{L^1(\Omega)}^{\frac{n-1}{n}}},
\end{equation*}
we have
\begin{equation}\label{HQ:30}
1\le\frac{-\zeta'_{n-1}(t)}{\omega_n^{1/n}}
\left(\int_0^{\omega_n\zeta_{n-1}(t)^n}f^*(s)\,ds\right)^{\frac1n}
-\frac{\omega_n^{\frac{n-1}{n}}}{4n\,\|f_0\|_{L^1(\Omega)}^{\frac{n-1}{n}}}\,
\frac{\mathcal D(t)^2}{(-\zeta'_{n-1}(t))^{n-1}} .
\end{equation}

\smallskip
\noindent\emph{Step 3: comparison with $v$.}
The first term on the right-hand side of \eqref{HQ:30} equals
$\frac{d}{dt}v(\zeta_{n-1}(t))$. Integrating \eqref{HQ:30} over $[0,\tau]$, and recalling
that $\zeta_{n-1}(0)=\zeta_{n-1}(\Omega)$ and $v(\zeta_{n-1}(\Omega))=0$, we get \eqref{HQ:31}. Taking $|x|=\zeta_{n-1}(\tau)$, so
that $\tau=u^\#_{n-1}(x)$, \eqref{HQ:31} gives
\[
\frac{\omega_n^{\frac{n-1}{n}}}{4n\,\|f_0\|_{L^1(\Omega)}^{\frac{n-1}{n}}}
\int_0^{u^\#_{n-1}(x)}\frac{\mathcal D(t)^2}{(-\zeta'_{n-1}(t))^{n-1}}\,dt\le v(x)-u^\#_{n-1}(x),
\]
and \eqref{HQ:HD} follows by taking the supremum over $\Omega^*_{n-1}$.
\end{proof}
We now apply the previous result, together with Lemma \ref{lem:geom}, to link the deficit $\mathcal{D}(t)$ with the Hausdorff asymmetry.
\begin{lemma}\label{HQ:thm}
Under the assumptions of Theorem \ref{Main Teo}, for every $0\le t<\|u\|_{L^\infty(\Omega)}$ such that $\zeta_{n-1}(t)\le\zeta_0(\Omega)$, we have
\begin{equation}\label{HQ:ratio}
\frac{\zeta_{n-1}(t)}{-\zeta'_{n-1}(t)}\ \ge\
\left(\frac{\|f_0\|_{L^1(\Omega)}}{|\Omega|}\right)^{\frac1n}\frac{\zeta_{n-1}(t)^2}{2}
-\|v-u^\#_{n-1}\|_{L^\infty(\Omega^*_{n-1})}.
\end{equation}
Moreover
\begin{equation}\label{HQ:final}
\int_0^{\|u\|_{L^\infty(\Omega)}}(-\zeta'_{n-1}(t))\,\mathcal D(t)\,dt\ \ge\
\gamma_n\,\zeta_{n-1}(\Omega)\,\alpha_{\mathcal H}(\Omega)^n,
\end{equation}
where the constant $\gamma_n$ is given in Theorem \ref{Main Teo}.
\end{lemma}
\begin{oss}
The assumption $\zeta_{n-1}(t)\le\zeta_0(\Omega)$ in the first
part of Lemma \ref{HQ:thm} is needed because \eqref{soglia}
holds only for $0\le\rho\le\zeta_0(\Omega)$.
Indeed, for $\rho>\zeta_0(\Omega)$, the left-hand side of
\eqref{soglia} is constant and equals
$\|f_0\|_{L^1(\Omega)}/\omega_n$, while the right-hand side
grows as $\rho^n$.
\end{oss}
\begin{proof}
\emph{Proof of \eqref{HQ:ratio}.} From \eqref{HQ:concave}-\eqref{normaTso} we have
\begin{equation*}
    \frac{\zeta_{n-1}(t)}{-\zeta_{n-1}'(t)}\geq \|u\|_{L^\infty(\Omega)}-t\geq v(0)-\|v-u^\#_{n-1}\|_{L^\infty(\Omega^*_{n-1})}-t.
\end{equation*}
We recall that \eqref{HQ:31} gives, in particular, $v(\zeta_{n-1}(t))\geq t$. Then the latter becomes
\begin{equation}\label{DisIntermedia}
    \begin{split}
        \frac{\zeta_{n-1}(t)}{-\zeta_{n-1}'(t)}&\geq v(0)-\|v-u^\#_{n-1}\|_{L^\infty(\Omega^*_{n-1})}-v(\zeta_{n-1}(t))\\
        &=\int_0^{\zeta_{n-1}(t)}\left(\frac{1}{\omega_n}\int_0^{\omega_n\rho^n}f^*(s)\,ds \right)^\frac{1}{n}\,d\rho-\|v-u^\#_{n-1}\|_{L^\infty(\Omega^*_{n-1})},
    \end{split}
\end{equation}
where the last equality follows from \eqref{v radiale}. If $\zeta_{n-1}(t)\leq\zeta_0(\Omega)$, then $\omega_n\rho^n\leq\abs{\Omega}$ for every $0\leq\rho\leq\zeta_{n-1}(t)$. Then, from Lemma \ref{f decrescente}, we have
\begin{equation}\label{soglia}
    \frac{1}{\omega_n}\int_0^{\omega_n\rho^n}f^*(s)\,ds\geq\frac{\rho^n\|f_0\|_{L^1(\Omega)}}{\abs{\Omega}}.
\end{equation}
Integrating the latter on $[0,\zeta_{n-1}(t)]$, we get
\begin{equation*}
    \int_0^{\zeta_{n-1}(t)}\left(\frac{1}{\omega_n}\int_0^{\omega_n\rho^n}f^*(s)\,ds\right)^\frac{1}{n}\,d\rho\geq\frac{\|f_0\|_{L^1(\Omega)}^\frac{1}{n}}{\abs{\Omega}^\frac{1}{n}}\frac{\zeta_{n-1}(t)^2}{2}.
\end{equation*}
Then, combining the latter with \eqref{DisIntermedia}, we find \eqref{HQ:ratio}.

\smallskip
\noindent\emph{Proof of \eqref{HQ:final}.}
By \eqref{HQ:supp1} and Lemma \ref{lem:geom}, it suffices to show that, for every
\(0<t<\|u\|_{L^\infty(\Omega)}\),
\begin{equation}\label{HQ:claim}
    \int_{\partial\Omega_t}
    H_{n-1}
    \left|\frac{1}{|\nabla u|}+\zeta'_{n-1}(t)\right|
    \,d\mathcal H^{n-1}
    \le
    n\omega_n\bigl(-\zeta'_{n-1}(t)\bigr)\mathcal D(t).
\end{equation}

Fix \(0<t<\|u\|_{L^\infty(\Omega)}\) and let
\[
    E_t:=
    \left\{
        x\in\partial\Omega_t:
        \frac{1}{|\nabla u(x)|}>-\zeta'_{n-1}(t)
    \right\}.
\]
By \eqref{HQ:identities} and the Gauss-Bonnet Theorem we have
\[
    \int_{\partial\Omega_t}
    H_{n-1}
    \left(\frac{1}{|\nabla u|}+\zeta'_{n-1}(t)\right)
    \,d\mathcal H^{n-1}=0.
\]
Since the expression in parentheses is positive exactly on \(E_t\), the latter gives
\begin{equation}\label{HQ:elem}
    \int_{\partial\Omega_t}
    H_{n-1}
    \left|\frac{1}{|\nabla u|}+\zeta'_{n-1}(t)\right|
    \,d\mathcal H^{n-1}
    =
    2\int_{E_t}
    H_{n-1}
    \left(\frac{1}{|\nabla u|}+\zeta'_{n-1}(t)\right)
    \,d\mathcal H^{n-1}.
\end{equation}

Now set
\[
    g_t(x):=
    \frac{1}{n\omega_n\bigl(-\zeta'_{n-1}(t)\bigr)|\nabla u(x)|}
    -
    \frac{|\nabla u(x)|^n}
    {n\displaystyle\int_{\Omega_t}\det(-D^2u)\,dx}.
\]
From \eqref{HQ:identities} and the definition \eqref{HQ:defD}, we have
\[
    \int_{\partial\Omega_t}H_{n-1}g_t\,d\mathcal H^{n-1}=0,
    \qquad
    \int_{\partial\Omega_t}H_{n-1}|g_t|\,d\mathcal H^{n-1}
    =\mathcal D(t).
\]
Then
\[
\begin{aligned}
    2\int_{E_t}H_{n-1}g_t\,d\mathcal H^{n-1}
    &=
    \int_{E_t}H_{n-1}g_t\,d\mathcal H^{n-1}
    -
    \int_{\partial\Omega_t\setminus E_t}
    H_{n-1}g_t\,d\mathcal H^{n-1}
    \\
    &\le
    \int_{\partial\Omega_t}H_{n-1}|g_t|\,d\mathcal H^{n-1}
    =\mathcal D(t).
\end{aligned}
\]

On the other hand, for every \(x\in E_t\), \eqref{TalentiIneq} gives
\[
    |\nabla u(x)|^n
    <
    \frac{1}{\bigl(-\zeta'_{n-1}(t)\bigr)^n}
    \le
    \frac{1}{\omega_n}
    \int_{\Omega_t}\det(-D^2u)\,dx.
\]
Therefore, on \(E_t\),
\[
\begin{aligned}
    \frac{1}{|\nabla u|}+\zeta'_{n-1}(t)
    &\le
    \frac{1}{|\nabla u|}
    -
    \frac{\omega_n\bigl(-\zeta'_{n-1}(t)\bigr)|\nabla u|^n}
    {\displaystyle\int_{\Omega_t}\det(-D^2u)\,dx}
    \\
    &=n\omega_n\bigl(-\zeta'_{n-1}(t)\bigr)g_t.
\end{aligned}
\]
Integrating the latter on \(E_t\) and using \eqref{HQ:elem}, we obtain
\[
\begin{aligned}
    \int_{\partial\Omega_t}
    H_{n-1}
    \left|\frac{1}{|\nabla u|}+\zeta'_{n-1}(t)\right|
    \,d\mathcal H^{n-1}
    &\le
    2n\omega_n\bigl(-\zeta'_{n-1}(t)\bigr)
    \int_{E_t}H_{n-1}g_t\,d\mathcal H^{n-1}
    \\
    &\le
    n\omega_n\bigl(-\zeta'_{n-1}(t)\bigr)\mathcal D(t),
\end{aligned}
\]
which is \eqref{HQ:claim}.

From \eqref{HQ:supp1} and \eqref{HQ:claim}, we get
\begin{equation}\label{HQ:supp}
    \frac{1}{n\omega_n}
    \int_{\mathbb S^{n-1}}
    \bigl|h_\Omega-\zeta_{n-1}(\Omega)\bigr|
    \,d\mathcal H^{n-1}
    \le
    \int_0^{\|u\|_{L^\infty(\Omega)}}
    \bigl(-\zeta'_{n-1}(t)\bigr)\mathcal D(t)\,dt.
\end{equation}
Combining \eqref{HQ:supp} and \eqref{Supp0}, we obtain \eqref{HQ:final}.
\end{proof}
We are now in position for the proof of the main result.
\begin{proof}[Proof of Theorem \ref{Main Teo}]
    We recall that \eqref{HQ:ratio} gives
    \begin{equation}\label{Inizio}
        \frac{\zeta_{n-1}(t)}{-\zeta_{n-1}'(t)}\geq \frac{\|f_0\|_{L^1(\Omega)}^\frac{1}{n}\zeta_{n-1}(t)^2}{2\abs{\Omega}^\frac{1}{n}}-\|v-u_{n-1}^\# \|_{L^\infty(\Omega^*_{n-1})}\geq\|v-u_{n-1}^\# \|_{L^\infty(\Omega^*_{n-1})}, 
    \end{equation}
    where the first inequality holds for any $t\in[0,\|u\|_{L^\infty(\Omega)})$ such that
    \begin{equation}\label{firstCond}
        \zeta_{n-1}(t)\leq\zeta_0(\Omega),
    \end{equation}
    and the second inequality holds for any $t\in[0,\|u\|_{L^\infty(\Omega)})$ such that
    \begin{equation}\label{secCond}
        \frac{\|f_0\|_{L^1(\Omega)}^\frac{1}{n}\zeta_{n-1}(t)^2}{2\abs{\Omega}^\frac{1}{n}}\geq2\|v-u_{n-1}^\# \|_{L^\infty(\Omega^*_{n-1})}.
    \end{equation}
    Now we want to find a condition for the existence of a value $t_0\in[0,\|u\|_{L^\infty(\Omega)})$ such that \eqref{firstCond}-\eqref{secCond} hold true.
    
    \smallskip
    \noindent\emph{Step 1: choice of $t_0$.} In this step we choose $t_0$ by minimizing the estimate below,
and find a condition ensuring that this choice satisfies
\eqref{firstCond}-\eqref{secCond}.\\
    Assume that a value $t_0\in[0,\|u\|_{L^\infty(\Omega)})$ with the required properties exists. Then \eqref{Inizio} gives
    \begin{equation*}
        \frac{\zeta_{n-1}(t_0)}{-\zeta_{n-1}'(t_0)}\geq\frac{\|f_0\|_{L^1(\Omega)}^\frac{1}{n}\zeta_{n-1}(t_0)^2}{4\abs{\Omega}^\frac{1}{n}}.
    \end{equation*}
    Since $\zeta_{n-1}(\,\cdot\,)$ is concave, we have that $-\zeta_{n-1}'(\,\cdot\,)$ is nondecreasing. Then the latter gives
    \begin{equation}\label{(1)}
        -\zeta_{n-1}'(t)\leq-\zeta_{n-1}'(t_0)\leq\frac{4\abs{\Omega}^\frac{1}{n}}{\|f_0\|_{L^1(\Omega)}^\frac{1}{n}\zeta_{n-1}(t_0)}\quad\text{for all }t\in(0,t_0]. 
    \end{equation}
    Since $0\le\mathcal D(t)\leq2$, we have
    \begin{equation}\label{dopot_0}
        \int_{t_0}^{\|u\|_{L^\infty(\Omega)}}-\zeta_{n-1}'(t)\mathcal D(t)\,dt\leq 2\zeta_{n-1}(t_0).
    \end{equation}
    Then using the Cauchy-Schwarz inequality and the monotonicity of $-\zeta'_{n-1}$ we have
    \begin{equation*}
        \begin{split}
            \int_0^{t_0}-\zeta_{n-1}'(t)\mathcal D(t)\,dt&\leq\left(\int_0^{t_0}\frac{\mathcal D(t)^2}{(-\zeta_{n-1}'(t))^{n-1}}\,dt\right)^\frac{1}{2}\left(\int_0^{t_0}(-\zeta_{n-1}'(t))^{n+1}\,dt\right)^\frac{1}{2}\\
            &\leq (-\zeta_{n-1}'(t_0))^\frac{n}{2}\left(\int_0^{t_0}\frac{\mathcal D(t)^2}{(-\zeta_{n-1}'(t))^{n-1}}\,dt\right)^\frac{1}{2}\left(\int_0^{t_0}-\zeta_{n-1}'(t)\,dt\right)^\frac{1}{2}.
        \end{split}
    \end{equation*}
    The third factor in the previous expression is bounded by the integral on the whole interval $(0,\|u\|_{L^\infty(\Omega)})$, that gives $\zeta_{n-1}(\Omega)^\frac{1}{2}$. Applying \eqref{(1)} to the first factor and \eqref{HQ:HD} to the second, we obtain
    \begin{equation}\label{prima t_0}
        \int_0^{t_0}-\zeta_{n-1}'(t)\mathcal D(t)\,dt\leq \frac{2^{n+1}n^\frac{1}{2}}{\omega_n^\frac{n-1}{2n}}\frac{\zeta_{n-1}(\Omega)^\frac{1}{2}\abs{\Omega}^\frac{1}{2}}{\|f_0\|_{L^1(\Omega)}^\frac{1}{2n}\zeta_{n-1}(t_0)^\frac{n}{2}}\|v-u_{n-1}^\# \|_{L^\infty(\Omega^*_{n-1})}^\frac{1}{2}.
    \end{equation}
    Combining \eqref{dopot_0}-\eqref{prima t_0} we get
    \begin{equation}\label{(2)}
        \int_0^{\|u\|_{L^\infty(\Omega)}}-\zeta_{n-1}'(t)\mathcal D(t)\,dt\leq\frac{2^{n+1}n^\frac{1}{2}}{\omega_n^\frac{n-1}{2n}}\frac{\zeta_{n-1}(\Omega)^\frac{1}{2}\abs{\Omega}^\frac{1}{2}}{\|f_0\|_{L^1(\Omega)}^\frac{1}{2n}\zeta_{n-1}(t_0)^\frac{n}{2}}\|v-u_{n-1}^\# \|_{L^\infty(\Omega^*_{n-1})}^\frac{1}{2}+2\zeta_{n-1}(t_0).
    \end{equation}
    We remark that in \eqref{(2)} the dependence on $t_0$ is only through $\zeta_{n-1}(t_0)$. Minimizing the right hand side of \eqref{(2)} with respect to $\zeta_{n-1}(t_0)$ gives
    \begin{equation}\label{(3)}
        \zeta_{n-1}(t_0)^\frac{n+2}{2}=\frac{2^{n-1}n^\frac{3}{2}}{\omega_n^\frac{n-1}{2n}}\frac{\zeta_{n-1}(\Omega)^\frac{1}{2}\abs{\Omega}^\frac{1}{2}}{\|f_0\|_{L^1(\Omega)}^\frac{1}{2n}}\|v-u_{n-1}^\# \|_{L^\infty(\Omega_{n-1}^*)}^\frac{1}{2}.
    \end{equation}
    With this value of $\zeta_{n-1}(t_0)$, \eqref{(2)} becomes
    \begin{equation}\label{(4)}
        \int_0^{\|u\|_{L^\infty(\Omega)}}-\zeta_{n-1}'(t)\mathcal D(t)\,dt\leq\frac{2(n+2)}{n}\zeta_{n-1}(t_0).
    \end{equation}
    Since $\zeta_{n-1}:[0,\|u\|_{L^\infty(\Omega)}]\rightarrow [0,\zeta_{n-1}(\Omega)]$ is a bijection, a number $t_0\in[0,\|u\|_{L^\infty(\Omega)}]$ such that \eqref{(3)} holds exists if and only if the right hand side of \eqref{(3)} is in $[0,\zeta_{n-1}(\Omega)^\frac{n+2}{2}]$, that is, if and only if
    \begin{equation*}
        \frac{2^{n-1}n^\frac{3}{2}}{\omega_n^\frac{n-1}{2n}}\frac{\zeta_{n-1}(\Omega)^\frac{1}{2}\abs{\Omega}^\frac{1}{2}}{\|f_0\|_{L^1(\Omega)}^\frac{1}{2n}}\|v-u_{n-1}^\# \|_{L^\infty(\Omega_{n-1}^*)}^\frac{1}{2}\leq\zeta_{n-1}(\Omega)^\frac{n+2}{2}.
    \end{equation*}
    The latter is equivalent to the following first condition
    \begin{equation}\tag{I}\label{I request}
        \|v-u_{n-1}^\# \|_{L^\infty(\Omega^*_{n-1})}\leq\frac{\omega_n^\frac{n-1}{n}}{4^{n-1}n^3}\frac{\zeta_{n-1}(\Omega)^{n+1}\|f_0\|_{L^1(\Omega)}^\frac{1}{n}}{\abs{\Omega}}.
    \end{equation}
    The second condition for the existence of $t_0$ comes from \eqref{firstCond} and gives
    \begin{equation*}
        \frac{2^{n-1}n^\frac{3}{2}}{\omega_n^\frac{n-1}{2n}}\frac{\zeta_{n-1}(\Omega)^\frac{1}{2}\abs{\Omega}^\frac{1}{2}}{\|f_0\|_{L^1(\Omega)}^\frac{1}{2n}}\|v-u_{n-1}^\# \|_{L^\infty(\Omega_{n-1}^*)}^\frac{1}{2}=\zeta_{n-1}(t_0)^\frac{n+2}{2}\leq \zeta_0(\Omega)^\frac{n+2}{2}=\frac{\abs{\Omega}^\frac{n+2}{2n}}{\omega_n^\frac{n+2}{2n}},
    \end{equation*}
    where the first equality is \eqref{(3)}. The latter is equivalent to the following second condition
    \begin{equation}\tag{II}\label{II request}
        \|v-u_{n-1}^\# \|_{L^\infty(\Omega^*_{n-1})}\leq \frac{1}{4^{n-1}n^3\omega_n^\frac{3}{n}}\frac{\|f_0\|_{L^1(\Omega)}^\frac{1}{n}\abs{\Omega}^\frac{2}{n}}{\zeta_{n-1}(\Omega)}.
    \end{equation}
    We remark that from the Alexandrov-Fenchel inequality we have
    \begin{equation*}
        \abs{\Omega}^\frac{n+2}{2n}\omega_n^{-\frac{n+2}{2n}}=\zeta_0(\Omega)^\frac{n+2}{2}\leq \zeta_{n-1}(\Omega)^\frac{n+2}{2}.
    \end{equation*}
    Combining \eqref{II request} with the latter, we find \eqref{I request}. Then the only condition so far is \eqref{II request}. The third condition comes from \eqref{secCond}, that is
    \begin{equation*}
        \|v-u_{n-1}^\# \|_{L^\infty(\Omega^*_{n-1})}\leq\frac{\|f_0\|_{L^1(\Omega)}^\frac{1}{n}\zeta_{n-1}(t_0)^2}{4\abs{\Omega}^\frac{1}{n}}.
    \end{equation*}
    Then substituting the expression for $\zeta_{n-1}(t_0)$ given by \eqref{(3)}, the latter is
    \begin{equation}\tag{III}\label{III request}
        \|v-u_{n-1}^\# \|_{L^\infty(\Omega^*_{n-1})}\leq \frac{2^\frac{2n-8}{n}n^\frac{6}{n}}{\omega_n^\frac{2(n-1)}{n^2}}\|f_0\|_{L^1(\Omega)}^\frac{1}{n}\abs{\Omega}^\frac{n-2}{n^2}\zeta_{n-1}(\Omega)^\frac{2}{n}.
    \end{equation}
    Then, starting from \eqref{II request} and applying the Alexandrov-Fenchel inequality, we have
    \begin{equation*}
        \begin{split}
            \|v-u_{n-1}^\# \|_{L^\infty(\Omega^*_{n-1})}&\leq\frac{2^\frac{2n-8}{n}n^\frac{6}{n}}{\omega_n^\frac{2(n-1)}{n^2}}\|f_0\|_{L^1(\Omega)}^\frac{1}{n}\abs{\Omega}^\frac{n-2}{n^2}\zeta_{n-1}(\Omega)^\frac{2}{n} n^{-\frac{3n+6}{n}}2^{-\frac{2n^2-8}{n}}\\
            &\leq \frac{2^\frac{2n-8}{n}n^\frac{6}{n}}{\omega_n^\frac{2(n-1)}{n^2}}\|f_0\|_{L^1(\Omega)}^\frac{1}{n}\abs{\Omega}^\frac{n-2}{n^2}\zeta_{n-1}(\Omega)^\frac{2}{n}.
        \end{split}
    \end{equation*}
    Then \eqref{II request} implies \eqref{III request}. Therefore, if \eqref{II request} holds, the level $t_0$ selected by \eqref{(3)}
exists and satisfies \eqref{firstCond}--\eqref{secCond}.

\smallskip
    \noindent\emph{Step 2: \eqref{II request} holds true.} Assume that \eqref{II request} holds true. Then $t_0$ exists and we can apply \eqref{(4)}
    \begin{equation*}
        \frac{2(n+2)}{n}\zeta_{n-1}(t_0)\geq\int_0^{\|u\|_{L^\infty(\Omega)}}-\zeta_{n-1}'(t)\mathcal D(t)\,dt\geq\gamma_n\zeta_{n-1}(\Omega)\alpha_\mathcal{H}(\Omega)^n,
    \end{equation*}
    where the last inequality is \eqref{HQ:final}. Replacing $\zeta_{n-1}(t_0)$ from \eqref{(3)}, the latter gives
    \begin{equation*}
        \|v-u_{n-1}^\# \|_{L^\infty(\Omega^*_{n-1})}\geq\frac{\gamma_n^{n+2}n^{n-1}\omega_n^\frac{n-1}{n}}{2^{3n}(n+2)^{n+2}}\frac{\|f_0\|_{L^1(\Omega)}^\frac{1}{n}\zeta_{n-1}(\Omega)^{n+1}}{\abs{\Omega}}\alpha_\mathcal{H}(\Omega)^{n(n+2)}.
    \end{equation*}
    Then combining the latter with \eqref{boundNormaInf di u} we find
    \begin{equation}\label{QuantitativaPrimoCaso}
        \begin{split}
            \frac{\|v-u_{n-1}^\# \|_{L^\infty(\Omega^*_{n-1})}}{\|u\|_{L^\infty(\Omega)}}&\geq \frac{\gamma_n^{n+2}n^{n-2}}{2^{3n+1}(n+2)^{n+2}}\frac{\zeta_{n-1}(\Omega)^{n+1}}{\zeta_0(\Omega)^{n+1}}\alpha_\mathcal{H}(\Omega)^{n(n+2)}\\&\geq \frac{\gamma_n^{n+2}n^{n-2}}{2^{3n+1}(n+2)^{n+2}}\alpha_\mathcal{H}(\Omega)^{n(n+2)},
        \end{split}
    \end{equation}
    where in the last inequality we have used $\zeta_0(\Omega)\leq\zeta_{n-1}(\Omega)$.

    \smallskip
    \noindent\emph{Step 3: \eqref{II request} does not hold and $\zeta_{n-1}(\Omega)\geq(4n+1)\zeta_0(\Omega)$.} By \eqref{normaTso} and the definition of $v$, we have that
    \begin{equation*}
        \begin{split}
            \|v-u_{n-1}^\# &\|_{L^\infty(\Omega^*_{n-1})}=v(0)-\|u\|_{L^\infty(\Omega)}\\
            &=\frac{1}{\omega_n^\frac{1}{n}}\int_0^{\zeta_{n-1}(\Omega)}\left(\int_0^{\omega_n\rho^n}f^*(s)\,ds \right)^\frac{1}{n}\,d\rho-\|u\|_{L^\infty(\Omega)}\\
            &=\frac{1}{\omega_n^\frac{1}{n}}\left[\int_0^{\zeta_0(\Omega)}\left(\int_0^{\omega_n\rho^n}f^*(s)\,ds \right)^\frac{1}{n}\,d\rho+\int_{\zeta_0(\Omega)}^{\zeta_{n-1}(\Omega)}\left(\int_0^{\omega_n\rho^n}f^*(s)\,ds \right)^\frac{1}{n}\,d\rho \right]-\|u\|_{L^\infty(\Omega)}.
        \end{split}
    \end{equation*}
    We remark that the inner integral in $ds$ is constantly equal to $\|f_0\|_{L^1(\Omega)}$ for every $\rho\in[\zeta_0(\Omega),\zeta_{n-1}(\Omega)]$. Moreover, on the first integral in $d\rho$ in the previous expression we can apply Lemma \ref{f decrescente}. Then the latter gives
    \begin{equation*}
        \|v-u_{n-1}^\# \|_{L^\infty(\Omega^*_{n-1})}\geq\frac{\|f_0\|_{L^1(\Omega)}^\frac{1}{n}}{\omega_n^\frac{1}{n}}\left(\zeta_{n-1}(\Omega)-\frac{\zeta_0(\Omega)}{2}\right)-\|u\|_{L^\infty(\Omega)}.
    \end{equation*}
    Applying \eqref{boundNormaInf di u} and the assumption of this step, from the latter we have 
    \begin{equation*}
        \|v-u_{n-1}^\# \|_{L^\infty(\Omega^*_{n-1})}\geq\frac{\|f_0\|_{L^1(\Omega)}^\frac{1}{n}}{\omega_n^\frac{1}{n}}\left(\zeta_{n-1}(\Omega)-\frac{4n+1}{2}\zeta_0(\Omega)\right)\geq \frac{\|f_0\|_{L^1(\Omega)}^\frac{1}{n}}{2\omega_n^\frac{1}{n}}\zeta_{n-1}(\Omega).
    \end{equation*}
    Combining the latter with \eqref{boundNormaInf di u} and using the assumption of this step, we get
    \begin{equation}\label{QuantitativaSecCaso}
        \begin{split}
            \frac{\|v-u_{n-1}^\# \|_{L^\infty(\Omega^*_{n-1})}}{\|u\|_{L^\infty(\Omega)}}&\geq\frac{\|f_0\|_{L^1(\Omega)}^\frac{1}{n}}{2\omega_n^\frac{1}{n}}\frac{\zeta_{n-1}(\Omega)}{\|u\|_{L^\infty(\Omega)}}\geq\frac{1}{4n}\frac{\zeta_{n-1}(\Omega)}{\zeta_0(\Omega)}\geq\frac{4n+1}{4n}\\&\geq\frac{(4n+1)\omega_{n-1}^{n(n+2)}}{4n^{1+n(n+2)}\omega_n^{n(n+2)}}\alpha_\mathcal{H}(\Omega)^{n(n+2)},
        \end{split}
    \end{equation}
    where the last inequality follows from \eqref{BoundHaus}.

    \smallskip
    \noindent\emph{Step 4: \eqref{II request} does not hold and $\zeta_{n-1}(\Omega)<(4n+1)\zeta_0(\Omega)$.}
    The negation of the condition \eqref{II request} together with \eqref{boundNormaInf di u} gives
    \begin{equation*}
        \frac{\|v-u_{n-1}^\# \|_{L^\infty(\Omega^*_{n-1})}}{\|u\|_{L^\infty(\Omega)}}>\frac{\|f_0\|_{L^1(\Omega)}^\frac{1}{n}\abs{\Omega}^\frac{2}{n}}{4^{n-1}n^3\omega_n^\frac{3}{n}\zeta_{n-1}(\Omega)}\frac{1}{\|u\|_{L^\infty(\Omega)}}\geq\frac{1}{2^{2n-1}n^4}\frac{\zeta_0(\Omega)}{\zeta_{n-1}(\Omega)}.
    \end{equation*}
    The latter together with the assumption of this step gives
    \begin{equation}\label{QuantitativaTerzoCaso}
        \begin{split}
            \frac{\|v-u_{n-1}^\# \|_{L^\infty(\Omega^*_{n-1})}}{\|u\|_{L^\infty(\Omega)}}&>\frac{1}{2^{2n-1}n^4}\frac{\zeta_0(\Omega)}{\zeta_{n-1}(\Omega)}>\frac{1}{2^{2n-1}n^4(4n+1)}\\&\geq\frac{\omega_{n-1}^{n(n+2)}}{2^{2n-1}n^{4+n(n+2)}(4n+1)\omega_n^{n(n+2)}}\alpha_\mathcal{H}(\Omega)^{n(n+2)},
        \end{split}
    \end{equation}
    where the last inequality follows from \eqref{BoundHaus}.

    \smallskip
    \noindent\emph{Step 5: Proof of \eqref{RisFinale}.} From \eqref{QuantitativaPrimoCaso}-\eqref{QuantitativaSecCaso}-\eqref{QuantitativaTerzoCaso} we have
    \begin{equation*}
        \begin{split}
            \frac{\|v-u_{n-1}^\# \|_{L^\infty(\Omega^*_{n-1})}}{\|u\|_{L^\infty(\Omega)}}&\geq\min\left\{\frac{\gamma_n^{n+2}n^{n-2}}{2^{3n+1}(n+2)^{n+2}},\,\frac{(4n+1)\omega_{n-1}^{n(n+2)}}{4n^{1+n(n+2)}\omega_n^{n(n+2)}},\, \frac{\omega_{n-1}^{n(n+2)}\omega_n^{-n(n+2)}}{2^{2n-1}n^{4+n(n+2)}(4n+1)}\right\}\alpha_\mathcal{H}(\Omega)^{n(n+2)}\\
            &= \frac{\gamma_n^{n+2}n^{n-2}}{2^{3n+1}(n+2)^{n+2}}\alpha_\mathcal{H}(\Omega)^{n(n+2)}.
        \end{split}
    \end{equation*}
\end{proof}
We now prove Corollary \ref{QuantABP}.
\begin{proof}[Proof of Corollary \ref{QuantABP}]
By \eqref{normaTso},
\[
\|v-u^\# _{n-1}\|_{L^\infty(\Omega^*_{n-1})}=v(0)-\|u\|_{L^\infty(\Omega)} .
\]
Moreover, for $0\le\rho\le\zeta_{n-1}(\Omega)$ we
have $\int_0^{\omega_n\rho^n}f^*(s)\,ds\le\|f_0\|_{L^1(\Omega)}$, and therefore
\[
v(0)=\int_0^{\zeta_{n-1}(\Omega)}\left(\frac1{\omega_n}\int_0^{\omega_n\rho^n}f^*(s)\,ds\right)^{\frac1n}d\rho
\le\zeta_{n-1}(\Omega)\left(\frac{\|f_0\|_{L^1(\Omega)}}{\omega_n}\right)^{\frac1n}
\le\frac{\mathrm{diam}(\Omega)}{2\,\omega_n^{1/n}}\,\|f_0\|_{L^1(\Omega)}^{\frac1n},
\]
where in the last inequality we used $2\zeta_{n-1}(\Omega)=w(\Omega)\le\mathrm{diam}(\Omega)$.
Hence
\[
\frac{\dfrac{\mathrm{diam}(\Omega)}{2\,\omega_n^{1/n}}\,\|f_0\|_{L^1(\Omega)}^{\frac1n}-\|u\|_{L^\infty(\Omega)}}{\|u\|_{L^\infty(\Omega)}}
\ \ge\ \frac{v(0)-\|u\|_{L^\infty(\Omega)}}{\|u\|_{L^\infty(\Omega)}}
=\frac{\|v-u^\# _{n-1}\|_{L^\infty(\Omega^*_{n-1})}}{\|u\|_{L^\infty(\Omega)}},
\]
and \eqref{RisABP} follows from \eqref{RisFinale}.
\end{proof}

\paragraph{Funding}
  The authors were partially supported by Gruppo Nazionale per l’Analisi Matematica, la Probabilità e le loro Applicazioni
(GNAMPA) of Istituto Nazionale di Alta Matematica (INdAM).\\
Francesco Salerno is partially supported by INdAM GNAMPA 2026 Project ``Processi di diffusione non-lineari: regolarità e classificazione delle soluzioni'', CUP E53C25002010001.

\paragraph{AI Usage Declaration}
All the results and proofs in this paper were obtained by the authors. Claude (Anthropic) was used only as an assistant for language editing and typesetting of the manuscript, and for checking computations. The authors take full responsibility for the content of the paper.
\paragraph{Conflict of interest}
   The authors declare that they have no conflicts of interest with respect to this work.
\paragraph{Data Availability} There are no data associated with this manuscript.

\bibliographystyle{plain}
\bibliography{biblio}

\Addresses

\end{document}